\documentclass[amsmath,amssymb,wasysym]{amsart}

 \newtheorem{thm}{Theorem}[subsection]
 
 \newtheorem{lem}[thm]{Lemma}
 
 \theoremstyle{definition}
 \newtheorem{defn}[thm]{Definition}
 \theoremstyle{remark}
 \newtheorem{rem}[thm]{Remark}
 \numberwithin{equation}{subsection}

\begin{document}

\title[Reidemeister Torsion of a Symplectic Complex and Pfaffian]
 {Reidemeister Torsion of a Symplectic Complex}

\author{Ya\c{s}ar S\"{o}zen}

\address{Rowan University,
 Department of Mathematics, Glassboro NJ 08028}

\email{sozen@rowan.edu}

\subjclass{15A90 ; 57M99}

\keywords{Reidemeister torsion, Symplectic Complex}

\date{December 4, 2004}

\dedicatory{}

\commby{}


\begin{abstract}
We consider a claim mentioned in \cite{Witten} pp 187 about the
relation between a symplectic chain complex with $\omega-$compatible
bases and Reidemeister Torsion of it. This is an explanation of it.
\end{abstract}

\maketitle

\section*{Introduction}
Even though, we approach Reidemeister torsion as a linear algebraic
object, it actually is a combinatorial invariant for the space of
representations of a compact surface into a fixed gauge group
\cite{Witten},\cite{Porti}.\\

More precisely, let $S$ be a compact surface with genus at least $2$
and without boundary, $G$ be a gauge group with its (semi-simple)
Lie algebra $\mathfrak{g}.$ Then, for a representation
$\rho:\pi_1(S)\to G,$ we can associate the corresponding adjoint
bundle $\left(
\begin{array}{c}
  \widetilde{S}\times_{\rho}\; \mathfrak{g} \\
  \downarrow \\
  S \\
\end{array}%
\right) $ over $S,$ i.e. $\widetilde{S}\times_{\rho}\;
\mathfrak{g}=\widetilde{S}\times\; \mathfrak{g}\;/\sim,$ where
$(x,t)$ is identified with all the elements in its orbit i.e.
$(\gamma\bullet x,\gamma\bullet t)$ for all $\gamma \in \pi_1(S),$
and where in the first component the element $\gamma\in \pi_1(S)$ of
the fundamental group of $S$ acts as a deck transformation, and
in the second component by conjugation by $\rho(\gamma).$\\

Suppose $K$ is a cell-decomposition of $S$ so that the adjoint
bundle $\widetilde{S}\times_{\rho}\; \mathfrak{g}$ on $S$ is trivial
over each cell. Let $\widetilde{K}$ be the lift of $K$ to the
universal covering $\widetilde{S}$ of $S.$ With the action of
$\pi_1(S)$ on $\widetilde{S}$ as deck transformation,
$C_{\ast}(\widetilde{K};\mathbb{Z})$ can be considered a
left$-\mathbb{Z}[\pi_1(S)]$ module and with the action of $\pi_1(S)$
on $\mathfrak{g}$ by adjoint representation, $\mathfrak{g}$ can be
considered as a left$-\mathbb{Z}[\pi_1(S)]$ module, where
$\mathbb{Z}[\pi_1(S)]$ is the integral group ring
$\left\{\displaystyle\sum_{i=1}^p m_i\gamma_i\; ; m_i\in
\mathbb{Z},\; \gamma_i\in \pi_1(S),\; p\in \mathbb{N}\right\}.$
\\

More explicitly, if $\displaystyle\sum_{i=1}^p m_i\gamma_i$ is in
$\mathbb{Z}[\pi_1(S)],$ $t$ is in $\mathfrak{g},$ and
$\displaystyle\sum_{j=1}^q n_j\;\sigma_j\in
C_{\ast}(\widetilde{S};\mathbb{Z}),$ then $
\left(\displaystyle\sum_{i=1}^p
m_i\gamma_i\right)\bullet\left(\displaystyle\sum_{j=1}^q
n_j\sigma_j\right)\stackrel{\text{defn}}{=}\displaystyle
\sum_{i,j}\;n_jm_i\;(\gamma_i\bullet\sigma_j ),$ where $\gamma_i$
acts on $\sigma_j\subset\widetilde{S}$ by deck transformation, and
$\left(\displaystyle\sum_{j=1}^q m_j\gamma_j\right)\bullet
t\stackrel{\text{defn}}{=}\displaystyle\sum_{j=1}^q
m_j\;(\gamma_j\bullet t),$ where $\gamma_j\bullet
t=Ad_{\rho(\gamma_j)}(t)=\rho(\gamma_j)t\rho(\gamma_j)^{-1}. $ \\

To talk about the tensor product $\displaystyle
C_{\ast}(\widetilde{K};\mathbb{Z} )\displaystyle\otimes
\mathfrak{g},$ we should consider the left
$\mathbb{Z}[\pi_1(S)]-$module $C_{\ast}(\widetilde{K};\mathbb{Z} )$
as a right $\mathbb{Z}[\pi_1(S)]-$module as
$\sigma\bullet\gamma\stackrel{defn}{=}\gamma^{-1}\bullet \sigma,$
where the action of $\gamma^{-1}$ is as a deck-transformation. Note
that the relation $\sigma\bullet\gamma\otimes
t=\sigma\otimes\gamma\bullet t$ becomes
$\gamma^{-1}\bullet\sigma\otimes t=\sigma\otimes \gamma\bullet t,$
equivalently $\sigma'\otimes t=\gamma\bullet\sigma'\otimes
\gamma\bullet t,$ where $\sigma'$ is $\gamma^{-1}\bullet \sigma.$ We
may conclude that tensoring with $\mathbb{Z}[\pi_1(S)]$ has the same
effect as factoring with $\pi_1(S).$ Thus,
 $C_{\ast}(K;Ad_{\rho})\stackrel{\text{defn} }{=}
C_{\ast}(\widetilde{K};\mathbb{Z} )\displaystyle\otimes_{\rho}
\mathfrak{g}$ is defined as the quotient $\displaystyle
C_{\ast}(\widetilde{K};\mathbb{Z})\displaystyle\otimes
\mathfrak{g}\;/\sim ,$ where the elements of the orbit
$\{\gamma\bullet\sigma\otimes\gamma\bullet t;\; \text{for\ all\ }\;
\gamma\in \pi_1(S)\}$ of $\sigma\otimes t$ are identified.\\

In this way, we obtain the following complex:
$$0\to C_2(K;Ad_{\rho})\stackrel{\partial_2\otimes \text{id}}{\longrightarrow }
C_1(K;Ad_{\rho})\stackrel{\partial_1\otimes
\text{id}}{\longrightarrow} C_0(K;Ad_{\rho})\to 0,
$$
where $\partial_i$ is the usual boundary operator. For this complex,
we can associate the homologies $H_{\ast}(K;Ad_{\rho}).$ Similarly,
the twisted cochains $C^{\ast}(K;Ad_{\rho})$ will result the
cohomologies $H^\ast(K;Ad_{\rho}),$ where
$C^{\ast}(K;Ad_{\rho})\stackrel{\text{defn}}{=}
\text{Hom}_{\mathbb{Z}[\pi_1(S)]}(C_{\ast}(\widetilde{K};\mathbb{Z}),\mathfrak{g})$
is the set of $\mathbb{Z}[\pi_1(S)]-$module homomorphisms from
$C_{\ast}(\widetilde{K};\mathbb{Z})$ into
$\mathfrak{g}.$ For more information, we refer \cite{Porti},\cite{Sozen},\cite{Witten}.\\

If $\rho,\rho':\pi_1(S)\to G$ are conjugate, i.e.
$\rho'(\cdot)=A\rho(\cdot)A^{-1}$ for some $A\in G,$ then
$C_{\ast}(K;Ad_{\rho})$ and $C_{\ast}(K;Ad_{\rho'})$ are isomorphic.
Similarly, the twisted cochains $C^{\ast}(K;Ad_{\rho})$ and
$C^{\ast}(K;Ad_{\rho'})$ are isomorphic. Moreover, the homologies
$H_{\ast}(K;Ad_{\rho})$ are independent of the cell-decomposition.
For details, see \cite{Sozen},\cite{Witten},\cite{Porti}.\\

If $\{ e^i_1,\cdots,e^i_{m_i}\}$ is a basis for the
$C_i(K;\mathbb{Z}),$ then $c_i:=\{
\widetilde{e^i_1},\cdots,\widetilde{e^i_{m_i}}\}$ will be a
$\mathbb{Z}[\pi_1(S)]-$basis for $C_i(\widetilde{K};\mathbb{Z}),$
where $\widetilde{e^i_j}$ is a lift of $e^i_j.$ If we choose a basis
$\mathcal{A}$ of $\mathfrak{g},$ then $c_i\otimes_{\rho}
\mathcal{A}$ will be a $\mathbb{C}-$basis for $C_i(K;Ad_{\rho}),$
called a \emph{geometric} basis for $C_i(K;Ad_{\rho}).$ Recall that
$C_i(K;Ad_{\rho})=C_i(\widetilde{K};\mathbb{Z})\otimes_{\rho}\mathfrak{g},$
is defined as the quotient
$C_i(\widetilde{K};\mathbb{Z})\otimes\mathfrak{g}\;/\sim,$ where we
identify the orbit $\{\gamma\bullet\sigma\otimes \gamma\bullet t;
\gamma\in \pi_1(S)\}$ of $\sigma\otimes t,$ and where the action of
the fundamental group in the first slot by deck-transformations, and
in the second slot by the conjugation with $\rho(\cdot).$\\

In this set-up, one can also define
$\text{Tor}(C_{\ast}(K;Ad_{\rho}),\{c_i\otimes_{\rho}
\mathcal{A}\}_{i=0}^2,\{\mathfrak{h}_i\}_{i=0}^2)$ the
\emph{Reidemiester torsion} of  the triple $K,Ad_{\rho},$ and
$\{\mathfrak{h_i}\}_{i=0}^2,$ where $\mathfrak{h_i}$ is a
$\mathbb{C}-$basis for $H_i(K;Ad_{\rho}).$ Moreover, one can easily
prove that this definition does not depend on the lifts
$\widetilde{e^i_j},$ conjugacy class of $\rho,$ and
cell-decomposition $K$ of the surface $S.$ Details can be found in
\cite{Sozen},\cite{Porti},\cite{Witten}.\\

Let $K,K'$ be dual cell-decompositions of $S$ so that $\sigma\in K$
$\sigma'\in K'$ meet at most once and moreover the diameter of each
cell has diameter less than, say, half of the injectivity radius of
$S.$ If we denote
$C_{\ast}=C_{\ast}(K;Ad_{\rho}),C'_{\ast}=C_{\ast}(K';Ad_{\rho}),$
then by the invariance of torsion under subdivision,
$\text{Tor}(C_{\ast})=\text{Tor}(C'_{\ast}).$ Let $D_{\ast}$ denote
the complex $C_{\ast}\oplus C'_{\ast}.$ Then, easily we have the
short-exact sequence $$0\to C_{\ast}\to D_{\ast}=C_{\ast}\oplus
C'_{\ast}\to C'_{\ast}\to 0.$$ The complex $D_{\ast}=C_{\ast}\oplus
C'_{\ast}$ can also be considered as a symplectic complex. Moreover,
in the case of irreducible representation $\rho:\pi_1(S)\to G,$
torsion $\text{Tor}(C_{\ast})$ gives a two-form on
$H^1(S;Ad_{\rho}).$  See \cite{Witten},\cite{Sozen}.\\

In this article, we will consider Reidemeister torsion as a linear
algebraic object and try to rephrase a statement mentioned in
\cite{Witten}. \\

The main result of the article is as stated in  \cite{Witten} pp 187
``the torsion of a symplectic complex  $(C_{\ast},\omega)$ computed
using a compatible set of measures is ``trivial" in the sense that''

\begin{thm}\label{Main_Theorem}
For a general symplectic complex $C_{\ast},$ if $\mathfrak{c}_p,$
$\mathfrak{h}_p$ are bases for $C_p,$ $H_p,$ respectively, then
$$\text{Tor}(C_{\ast},\{\mathfrak{c}_p\}_{p=0}^n,\{\mathfrak{h}_p\}_{p=0}^n)=
\left(\displaystyle\prod_{p=0}^{\frac{n}{2}-1} (\det
[\omega_{p,n-p}])^{(-1)^p}\right)\cdot
   \left(\sqrt{\det[\omega_{\frac{n}{2},\frac{n}{2}}]}\;\right)^{(-1)^{\frac{n}{2}}},$$
where $\det[\omega_{p,n-p}]$ is the determinant of the matrix of the
non-degenerate pairing $[\omega_{p,n-p}]:H_p(C)\times H_{n-p}(C)\to
\mathbb{R}$ in bases $\mathfrak{h}_p,\mathfrak{h}_{n-p}.$\\
\end{thm}

For topological application of this, we refer \cite{Sozen},\cite{Witten}.\\

The plan of paper is as follows. In section \S 1, we will give the
definition of Reidemeister torsion for a general complex $C_{\ast}$
and recall some properties. See \cite{Milnor},\cite{Porti} for more
information. In section \S 2, we will explain torsion using Witten's
notation \cite{Witten}. Finally, symplectic complex will be
explained in section \S 3  and also the proof of main result Theorem~\ref{Main_Theorem}.\\

\section{Reidemeister Torsion of a general Chain Complex}

Let $C_{\ast }=(C_n\stackrel{\partial_n}{\longrightarrow}
C_{n-1}\longrightarrow\cdots\longrightarrow C_1
\stackrel{\partial_1}{\longrightarrow} C_0\longrightarrow 0)$ be a
chain complex of a finite dimensional vector spaces over
$\mathbb{R}$ or $\mathbb{C}.$ Let $H_p=Z_p/B_p$ denote the
homologies of the complex, where
$B_p=\text{Im}\{\partial_{p+1}:C_{p+1}\to C_p\},$
$Z_p=\text{ker}\{\partial_p:C_p\to C_{p-1}\},$ respectively.\\

If we start with bases ${\mathfrak b}_p=\{b_p^1,\cdots,b_p^{m_p}\}$
for $B_p$, and ${\mathfrak h}_p=\{h_p^1,\cdots,h_p^{n_p}\}$ for
$H_p$, a new basis for $C_p$ can be obtained by considering the
following short-exact sequences:

\begin{eqnarray}
 \label{Z_p-C_p-B_p-1} 0&\to &  Z_p  \hookrightarrow  C_p \rightarrow   B_{p-1} \to  0\\
 \label{B_p-Z_p-H_p} 0&\to &  B_p  \hookrightarrow  Z_p \rightarrow   H_p \to  0
\end{eqnarray}
where the first row is a result of $1^{\text{st}}$-Isomorphism
Theorem and the second follows simply from the definition of
$H_p.$\\

 Starting with \ref{B_p-Z_p-H_p} and a section
$\overline{s_p}:H_p\to Z_p,$ then $Z_p$ will have a basis
${\mathfrak b}_p\oplus \overline{s_p}(\mathfrak{h}_p).$ Using
\ref{Z_p-C_p-B_p-1} and a section $s_p:B_{p-1}\to C_p,$ $C_p$ will
have a basis $\mathfrak{b}_p\oplus
\overline{s_p}(\mathfrak{h}_p)\oplus s_p(\frak{b}_{p-1}).$\\

If $V$ is a vector space with bases $\mathfrak{e}$ and
$\mathfrak{f},$ then we will denote $[\mathfrak{f},\mathfrak{e}]$
for the determinant of the change-base-matrix
$T_{\mathfrak{e}}^{\mathfrak{f}}$ from
$\mathfrak{e}$ to $\mathfrak{f}.$ \\

\begin{defn}For $p=1,\cdots,n$, let $\mathfrak{c}_p,$ $\mathfrak{b}_p,$ and $\mathfrak{h}_p$ be bases
for $C_p,$ $B_p$ and $H_p,$ respectively.
$\text{Tor}(C_{\ast},\{\mathfrak{c}_p\}_{p=0}^n,\{\mathfrak{h}_p\}_{p=0}^n)=\displaystyle\prod_{p=0}^n
\left[\mathfrak{b}_p\oplus \overline{s_p}(\mathfrak{h}_p)\oplus
s_p(\frak{b}_{p-1}), \mathfrak{c}_p\right]^{(-1)^{(p+1)}}$ is called
the \emph{torsion of the complex} $C_{\ast}$ with respect to bases
$\{\mathfrak{c}_p\}_{p=0}^n,\{\mathfrak{h}_p\}_{p=0}^n,$ \\
\end{defn}

Milnor \cite{Milnor} proved that torsion does not depend on neither
the bases $\frak{b}_p,$ nor the sections $s_p,\overline{s_p}.$ In
other words, it is well-defined.

\begin{rem}\label{FactsAboutTorsion} If we choose another bases
    $\frak{c}'_p,\frak{h}'_p$ respectively for $C_p$ and $H_p$, then
    an easy computation shows that
$$\text{Tor}(C_{\ast},\{\mathfrak{c}'_p\}_{p=0}^n,\{\mathfrak{h}'_p\}_{p=0}^n)=
\displaystyle\prod_{p=0}^n\left(\dfrac{[\frak{c}'_p,\frak{c}_p]}{[\frak{h}'_p,\frak{h}_p]}\right)^{(-1)^p}\cdot
\text{Tor}(C_{\ast},\{\mathfrak{c}_p\}_{p=0}^n,\{\mathfrak{h}_p\}_{p=0}^n).$$
\end{rem}

This follows easily from the fact that torsion is independent of
$\frak{b}_p$ and sections $s_p,\overline{s_p}.$ For example, if
$[\frak{c}'_p,\frak{c}_p]=1,$ and $[\frak{h}'_p,\frak{h}_p]=1,$ then
they produce the same torsion.\\

If we have a short-exact sequence of chain complexes $ 0\to A_{\ast
}\stackrel{\imath}{\to} B_{\ast }\stackrel{\pi}{\to} D_{\ast }\to
0,$ then we also have a long-exact sequence of vector space
$C_{\ast}$
$$\cdots \to H_p(A)\stackrel{\imath_{\ast}}{\to }
H_p(B)\stackrel{\pi_{\ast}}{\to } H_p(D)\stackrel{\Delta}{\to}
H_{p-1}(A)\to\cdots $$ i.e. an acyclic (or exact) complex $C_{\ast}$
of length $3n+2$ with $C_{3p}=H_p(D_{\ast}),$
$C_{3p+1}=H_p(A_{\ast})$ and $C_{3p+2}=H_p(B_{\ast}).$ In
particular, the bases $\mathfrak{h}_p(D_{\ast}),$
$\mathfrak{h}_p(A_{\ast}),$ and $\mathfrak{h}_p(B_{\ast})$ will
serve as bases for $C_{3p},C_{3p+1},$ and $C_{3p+2},$ respectively.

\begin{thm}(Milnor \cite{Milnor})\label{Mil}
Using the above setup, let
$\mathfrak{c}^A_p,\mathfrak{c}^B_p,\mathfrak{c}^D_p$ be bases for
$A_p,B_p,D_p,$ respectively, and let
$\mathfrak{h}^A_p,\mathfrak{h}^B_p,\mathfrak{h}^D_p$ be bases for
the corresponding homologies $H_p(A),H_p(B),$ and $H_p(D).$ If,
moreover, the bases
$\mathfrak{c}^A_p,\mathfrak{c}^B_p,\mathfrak{c}^D_p$ are compatible
in the sense that $[\mathfrak{c}^B_p,\mathfrak{c}^A_p\oplus
\widetilde{\mathfrak{c}^D_p}]=\pm 1$ where
$\pi\left(\widetilde{\mathfrak{c}^D_p}\right)=\mathfrak{c}^D_p,$
then $\text
{Tor}(B_{\ast},\{\mathfrak{c}^B_p\}_{p=0}^n,\{\mathfrak{h}^B_p\}_{p=0}^n)$
$=\text{Tor}(A_{\ast},\{\mathfrak{c}^A_p\}_{p=0}^n,\{\mathfrak{h}^A_p\}_{p=0}^n)
\cdot
\text{Tor}(D_{\ast},\{\mathfrak{c}^D_p\}_{p=0}^n,\{\mathfrak{h}^D_p\}_{p=0}^n)\cdot
\text{Tor}(C_{\ast},\{\mathfrak{c}_{3p}\}_{p=0}^{3n+2},\{0\}_{p=0}^{3n+2}).$\\

 \hfill $\Box$
\end{thm}

\section{Reidemeister Torsion using Witten's notations}
Let $V$ be a vector space of dimension $k$ over $\mathbb{R}.$ Let
$\text{det}(V)$ denote the top exterior power $\bigwedge^kV$ of $V.$
A \emph{measure} on $V$ is a non-zero functional
$\alpha:\text{det}(V)\to\mathbb{R}$ on $\text{det}(V),$ i.e. $\alpha
\in {\text{det}(V)}^{-1},$ $-1$ denotes the dual space.\\

Recall that the isomorphism between ${\text{det}(V)}^{-1}$ and
$\text{det}(V^{\ast})$ is given by the pairing
$<\cdot,\cdot>:\text{det}(V^{\ast})\times
\text{det}(V)\to\mathbb{R},$ defined by
$$<f_1^{\ast}\wedge\cdots\wedge f_k^{\ast},e_1\wedge\cdots\wedge
e_k>=\text{det}\left[f_i^{\ast}(e_j)\right], $$
i.e.
$[\mathfrak{f},\mathfrak{e}]$ the determinant of the
change-base-matrix from basis $\mathfrak{e}=\{ e_1,\cdots,e_k\}$ to
$\mathfrak{f}=\{ f_1,\cdots,f_k\},$ where $f^{\ast}_i$  is the dual
element corresponding to $f_i,$ namely,
$f^{\ast}_i(f_j)=\delta_{ij}.$ Below $(v_1\wedge\cdots\wedge
v_k)^{-1}$ will denote $(v_1)^{\ast}\wedge\cdots\wedge (v_k)^{\ast}$\\

Note also that $<f^{\ast}_1\wedge\cdots\wedge f^{\ast}_k\;
,\;e_1\wedge\cdots e_k> =<e^{\ast}_1\wedge\cdots\wedge e^{\ast}_k\;
,\;f_1\wedge\cdots f_k>^{-1},$ i.e.
$[\mathfrak{f},\mathfrak{e}]=[\mathfrak{e},\mathfrak{f}]^{-1}.$ So,
using the pairing, $[\mathfrak{f},\;\bullet\;]\;$ can be considered
a linear functional on $\text{det}(V)$ and
$[\;\bullet\;,\mathfrak{e}]\;$ can be considered a linear functional
on $\text{det}(V^{\ast}).$\\

Let $C_{\ast}:0\to C_n\stackrel{\partial_n}{\to} C_{n-1}\to\cdots\to
C_1\stackrel{\partial_1}{\to} C_0\to 0$ be a chain complex of finite
dimensional vector spaces with \emph{volumes} $\alpha_p\in
{\text{det}(V)}^{-1},$ i.e.
$\alpha_p=(c^p_1)^{\ast}\wedge\cdots\wedge (c^p_{m_p})^{\ast}$ for
some basis $\{ c^p_1,\cdots, c^p_{m_p}\}$ for $C_p.$ If, moreover,
we assume that $C_{\ast}$ is exact (or acyclic), then
$H_p(C_{\ast})=0$ for all $p.$ In particular, we have the short
exact sequence
$$0\to\underbrace{\text{Im}\{\partial_{p+1}:C_{p+1}\to
C_p\}}_{B_p}\stackrel{i_p}{\hookrightarrow}
C_p\stackrel{\partial_p}{\to}\underbrace{\text{Im}\{\partial_p:C_p\to
C_{p-1}\}}_{B_{p-1}}\to 0.$$

Let $\{ b^p_1,\cdots,b^p_{k_p}\},$
$\{b^{p-1}_1,\cdots,b^{p-1}_{k_{p-1}}\}$ be bases for $B_p,B_{p-1},$
respectively. Then,
$\{b^p_1,\cdots,b^p_{k_p},\widetilde{b^{p-1}_1},\cdots,\widetilde{b^{p-1}_{k_{p-1}}}\}$
is a basis for $C_p,$ where
$\partial_p(\widetilde{b_{p-1}^i})=b_{p-1}^i$ and thus
$b^p_1\wedge\cdots\wedge b^p_{k_p}\wedge
\widetilde{b^{p-1}_1}\wedge\cdots\wedge\widetilde{b^{p-1}_{k_{p-1}}}$
is a basis for $\text{det}(C_p).$\\

 If $u$ denotes
$\displaystyle\bigotimes_{p=0}^{n}(b^p_1\wedge\cdots\wedge
b^p_{k_p}\wedge
\widetilde{b^{p-1}_1}\wedge\cdots\wedge\widetilde{b^{p-1}_{k_{p-1}}})^{(-1)^p},$
then $u$ is an element of $\displaystyle\bigotimes_{p=0}^{n}
(\text{det}(C_p))^{(-1)^p},$ where the exponent $(-1)$ denotes the
dual of the vector space. Then, E. Witten describes the torsion as:
\begin{eqnarray*}
  \text{Tor}(F_{\ast}) &=& < u,\displaystyle\bigotimes_{p=0}^{n}\;
 \alpha_p^{(-1)^p}> \\
 &=& \prod_{p=0}^n < b^p_1\wedge\cdots\wedge b^p_{k_p}\wedge
\widetilde{b^{p-1}_1}\wedge\cdots\wedge\widetilde{b^{p-1}_{k_{p-1}}},(c^p_1)^{\ast}\wedge\cdots\wedge
(c^p_{m_p})^{\ast}  >^{(-1)^p},
\end{eqnarray*}
which is nothing but
$\displaystyle\prod_{p=0}^n[\;\{c^p_1,\cdots,c^p_{m_p}\}\;,\{b^p_1,\cdots,
b^p_{k_p},
\widetilde{b^{p-1}_1},\cdots,\widetilde{b^{p-1}_{k_{p-1}}}\}\;
]^{(-1)^p}$
or\\
$\displaystyle\prod_{p=0}^n(\;[\;\{
b^p_1,\cdots,b^p_{k_p},\widetilde{b^{p-1}_1},
\cdots,\widetilde{b^{p-1}_{k_{p-1}}}\}\;,\;\{c^p_1,\cdots,c^p_{m_p}\}\;]^{(-1)})^{(-1)^p}.$
The last term coincides with the
$\text{Tor}(C_{\ast},\{c_p\}_{p=0}^n,\{0\}_{p=0}^n)$ defined in
previous section.\\

We will now explain how a general chain complex can be (unnaturally)
as a direct sum of an two chain complexes, one of which is exact and
the other is $\partial-$zero.

\begin{thm}\label{Unnatural_Splitting_Of_General_Complex}
If $C_{\ast}:0\to C_n\stackrel{\partial_n}{\to} C_{n-1}\to\cdots\to
C_1\stackrel{\partial_1}{\to} C_0\to 0$ is any chain complex, then
it can be splitted as $C_{\ast}=C'_{\ast}\oplus C''_{\ast},$ where
$C'_{\ast}$ is exact, and $C''_{\ast}$ is $\partial-$zero.
\end{thm}

\begin{proof} Consider the short-exact sequences
$$\begin{array}{cccccccc}
  0 & \to & \ker\partial_p & \hookrightarrow & C_p & \stackrel{\partial_p}{\rightarrow }  & \text{Im} \partial_p  & \to  0 \\
  0 & \to & \text{Im} \partial_{p+1} & \hookrightarrow & \ker\partial_p & \stackrel{\pi_p}{\rightarrow}  & H_p  & \to  0. \\
 \end{array}$$

If $\ell_p:\text{Im} \partial_p\to C_p,$ and $s_p:H_p\to
\ker\partial_p$ are sections, i.e.
$\partial_p\circ\ell_p=id_{\text{Im}_{\partial_p}},$ and $\pi_p\circ
s_p=id_{H_p(C)},$ then $C_p$ is equal to
$\ker\partial_p\oplus\ell_p(\text{Im}\partial_p)$ or
$\text{Im}\partial_{p+1}\oplus
s_p(H_p)\oplus\ell_p(\text{Im}\partial_p).$ Define
$C'_p:=\text{Im}\partial_{p+1}\oplus \ell_p(\text{Im} \partial_p)$
and $C''_p:=s_p(H_p).$ Restricting $\partial_p:C_p\to C_{p-1}$ to
these, we obtain two chain complexes $(C'_{\ast},\partial'_{\ast})$
$(C''_{\ast},\partial''_{\ast}).$\\

As $C''_p$ is a subspace of $\ker\partial_p,$ $\partial''_p:C''_p\to
C''_{p-1}$ is the zero map, i.e. $C''_{\ast}$ is $\partial-$zero
chain complex. Note also $\ker\{\partial''_p:C''_p\to C''_{p-1}\}$
equals to $C''_p$ and $\text{Im}\{\partial''_{p+1}:C''_{p+1}\to
C''_p\}$ is $\{0\}.$ Then, $H_p(C''_{\ast})=C''_p/\{ 0\}$ is
isomorphic to $H_p(C),$ because $C''_p=s_p(H_p(C))$ is isomorphic to
$H_p(C).$\\

The exactness of $(C'_{\ast},\partial'_{\ast})$ can be seen as
follows: Since $\text{Im}\partial_{p+1}$ is a subspace of
$\ker\partial_p,$ the image of $\text{Im}\partial_{p+1}$ under
$\partial'_p$ is zero. Hence, $\ker\{\partial'_p:C'_p\to C'_{p-1}\}$
equals to $\text{Im}\{\partial_{p+1}:C_{p+1}\to C_p\}.$ Since
$\partial_p\circ\ell_p=id_{\text{Im}_{\partial_p}},$ and
$\partial'_p:C'_p\to C'_{p-1}$ is the restriction of
$\partial_p:C_p\to C_{p-1},$ then $\text{Im}\{\partial'_p:C'_p\to
C'_{p-1}\}$ equals to $\text{Im}\{\partial_p:C_p\to C_{p-1}\}.$
Similarly, $\text{Im}\{\partial'_{p-1}:C'_{p-1}\to
C'_{p-2}\}=\text{Im}\{\partial_{p-1}:C_{p-1}\to C_{p-2}\}$ and
$\ker\{\partial'_{p-1}:C'_{p-1}\to
C'_{p-2}\}=\text{Im}\{\partial_{p}:C_{p}\to C_{p-1}\},$ because
$\text{Im}\partial_p$ is a subspace of $\ker\partial_{p-1}$ and
$\ell_{p-1}$ is a section of $\partial_{p-1}:C_{p-1}\to
\text{Im}\partial_{p-1}.$ Consequently,
$\text{Im}\{\partial'_p:C'_p\to
C'_{p-1}\}=\ker\{\partial'_{p-1}:C_{p-1}\to
C_{p-2}\}=\text{Im}\partial_p$ and we have the exactness of
$C'_{\ast}.$\\

This concludes Theorem~\ref{Unnatural_Splitting_Of_General_Complex}.\\
\end{proof}

In the next result, we will explain Witten's remark on
(\cite{Witten} pp 185) how torsion $\text{Tor}(C_{\ast})$ of a
general complex can be interpreted as an element of the dual of the
one dimensional vector space $\otimes_{p=0}^n\;
(\det(H_p(C)))^{(-1)^p}.$

\begin{thm}\label{Torsion_Is_In_The_Dual}
$\text{Tor}(C_{\ast})$ of a general complex is as an element of the
dual of the one dimensional vector space
$\displaystyle\bigotimes_{p=0}^{n}\;(\det(H_p(C)))^{(-1)^p}.$
\end{thm}

\begin{proof}
Let $C_{\ast}$ be a general chain complex of finite dimensional
vector spaces with volumes $\alpha_p\in(\text{det}C_p)^{-1},$ i.e.
$\alpha_p=(c_p^1)^{\ast}\wedge\cdots\wedge(c_p^{i_p})^{\ast},$ for
some basis $\mathfrak{c}_p=\{ c_p^1,\cdots,c_p^{i_p}\}$ of $C_p.$
Let $C_{\ast}=C'_{\ast}\oplus C''_{\ast}$ be the above unnatural
splitting of $C_{\ast}$ i.e. $C'_p=\text{Im}\partial_{p+1}\oplus
\ell_p(\text{Im}\partial_p)$ and $C''_p=s_p(H_p(C)),$ where
$\ell_p:\text{Im} \partial_p\to C_p$ is the lift of
$\partial_p:C_p\to \text{Im} \partial_p$ and $s_p:H_p\to
\ker\partial_p$ is the lift of $\pi_p:\ker\partial_p\to H_p(C)$
 used in Theorem~\ref{Unnatural_Splitting_Of_General_Complex}.\\

Since $C_p=\text{Im}\partial_{p+1}\oplus s_p(H_p)\oplus
\ell_p(\text{Im}\partial_p),$ we can break the basis
$\mathfrak{c}_p$ of $C_p$ into three blocks as
$\mathfrak{c}^p_1\sqcup \mathfrak{c}^p_2\sqcup \mathfrak{c}^p_3,$
where $\mathfrak{c}^p_1$ generates $\text{Im}\partial_{p+1},$
$\mathfrak{c}^p_2$ is basis for $s_p(H_p(C))$ i.e.
$[\mathfrak{c}^p_2]=\pi_p(\mathfrak{c}^p_2)$ generates $H_p(C),$ and
$\mathfrak{c}^p_3$ is a basis for $\text{Im}\partial_p.$ As the
determinant of change-base-matrix from $\mathfrak{c}_p$ to
$\mathfrak{c}_p$ is $1,$ the bases $\mathfrak{c}^p_2,$
$\mathfrak{c}_p=\mathfrak{c}^p_1\sqcup\mathfrak{c}^p_2\sqcup\mathfrak{c}^p_3,$
and $\mathfrak{c}^p_1\sqcup\mathfrak{c}^p_3$ for $C''_p,C_p,C'_p,$
will be compatible with the short-exact sequence of complexes

$$0\to
C''_{\ast}\hookrightarrow C_{\ast}=C''_{\ast}\oplus C'_{\ast}\to
C'_{\ast}\to 0,$$

\noindent where we consider the inclusion as section $C'_p\to C_p.$
Note also that $H_p(C'')=C''_p/0$ i.e. $s_p(H_p(C))$ which is
isomorphic to
$H_p(C).$\\

 By Milnor's result Theorem~\ref{Mil}, we have
$\text{Tor}(C_{\ast},\{\mathfrak{c}_p\}_{p=0}^n,\{\mathfrak{h}_p\}_{p=0}^n)$
is the product of
$\text{Tor}(C''_{\ast},\{\mathfrak{c}^2_p\}_{p=0}^n,\{s_p(\mathfrak{h}_p)\}_{p=0}^n),$
$\text{Tor}(C'_{\ast},\{\mathfrak{c}^1_p\sqcup
\mathfrak{c}^3_p\}_{p=0}^n,\{0\}_{p=0}^n),$
 and $\text{Tor}(\mathcal{H}_{\ast}),$
 where $\mathcal{H}_{\ast}$ is
 the long-exact sequence obtained from the above short-exact of
 chain complexes.\\

  More precisely, $\mathcal{H}_{\ast}:0\to H_n(C'')\to H_n(C)\to H_n(C')\to
H_{n-1}(C'')\to H_{n-1}(C)\to H_{n-1}(C')\to\cdots \to H_0(C'')\to
H_0(C)\to H_0(C')\to 0.$ Since $C'_{\ast}$ is exact, then
$\mathcal{H}_{\ast}$ is the long exact-sequence $0\to H_n(C'')\to
H_n(C)\to 0\to H_{n-1}(C'')\to H_{n-1}(C)\to 0\to\cdots \to 0 \to
H_0(C'')\to H_0(C)\to 0\to 0.$ Using the isomorphism $H_p(C)\to
H_p(C''),$ namely $s_p$  as section, we conclude that
$\text{Tor}(\mathcal{H}_{\ast},\{
s_p(\mathfrak{h}_p),\mathfrak{h}_p,0\}_{p=0}^{n},\{0\}_{p=0}^{3n+2})=1.$\\

Moreover, we can also verify that
$\text{Tor}(C'_{\ast},\{\mathfrak{c}^1_p\sqcup
\mathfrak{c}^3_p\}_{p=0}^n,\{0\}_{p=0}^n)=1$ as follows:\\
$$0\to \ker\{\partial'_p;C'_p\to C'_{p-1}\} \hookrightarrow C'_p
\stackrel{\partial'_p\stackrel{i.e.}{=}\partial_p}{\rightarrow}
\text{Im}\{\partial'_p:C'_p\to C'_{p-1}\} \to 0,$$ where
$\ker\{\partial'_p:C'_p\to C'_{p-1}\}$ is
$\text{Im}\{\partial_{p+1}:C_{p+1}\to C_p\}$ and $\text{Im}
\{\partial'_p:C'_p\to C'_{p-1}\}$ is $\text{Im}\{\partial_p:C_p\to
C_{p-1}\}.$ If we consider the section $\ell_p,$ then we also have
$\text{Tor}(C'_{\ast},\{\mathfrak{c}^1_p\sqcup
\mathfrak{c}^3_p\}_{p=0}^n,\{0\}_{p=0}^n)=1.$\\

Therefore,
$\text{Tor}(C_{\ast},\{\mathfrak{c}_p\}_{p=0}^n,\{\mathfrak{h}_p\}_{p=0}^n)$
is equal to
$\text{Tor}(C''_{\ast},\{\mathfrak{c}^2_p\}_{p=0}^n,\{s_p(\mathfrak{h}_p)\}_{p=0}^n)$
i.e.
$\displaystyle\prod_{p=0}^n[s_p(\mathfrak{h}_p),\mathfrak{c}^2_p]^{(-1)^p},$
where $[s_p(\mathfrak{h}_p),\mathfrak{c}^2_p]$ is the determinant of
the change-base-matrix from $\mathfrak{c}^2_p$ to
$s_p(\mathfrak{h}_p)$ of $C''_p=s_p(H_p(C)).$ Recall that
$s_p:H_p(C)\to \ker\partial_p$ is the lift of
$\pi_p:\ker\partial_p\to H_p(C).$ So, $[\mathfrak{c}^2_p],$ i.e.
$\pi_p(\mathfrak{c}_p),$ and $\mathfrak{h}_p=[s_p(\mathfrak{h}_p)]$
are bases for $H_p(C).$ Since $s_p$ is isomorphism onto its image,
change-base-matrix from $\mathfrak{c}^2_p$ to $s_p(\mathfrak{h}_p)$
coincides with the one from $[\mathfrak{c}^2_p]$ to
$\mathfrak{h}_p.$\\

As a result, we obtained that
\begin{eqnarray*}
  \text{Tor}(C_{\ast},\{\mathfrak{c}_p\}_{p=0}^n,\{\mathfrak{h}_p\}_{p=0}^n)
   &=& \displaystyle\prod_{p=0}^n[\mathfrak{h}_p,\;[\mathfrak{c}^2_p]\;]^{(-1)^p} \\
   &=& [\mathfrak{h}_0,\;[\mathfrak{c}^2_0]\;]\cdot
   [\mathfrak{h}_1,\;[\mathfrak{c}^2_1]\;]^{-1}\cdots
   [\mathfrak{h}_n,\;[\mathfrak{c}^2_n]\;]^{(-1)^n}.
\end{eqnarray*}

For $p$ even, $[\mathfrak{h}_p,\;[\mathfrak{c}^2_p]\;]^{(-1)^p}$ is
$[\mathfrak{h}_p,\;[\mathfrak{c}^2_p]\;],$ for $p$ odd, it is
$[\mathfrak{h}_p,\;[\mathfrak{c}^2_p]\;]^{-1}$ or $[\;
[\mathfrak{c}^2_p]\;,\;\mathfrak{h}_p\;].$\\

 By the remark at the
beginning of \S 2, for even $p$'s,
$[\;\bullet\;,\;[\mathfrak{c}^2_p]\;]$ is linear functional on
$\det(H_p(C)),$ and for odd $p$'s,
$[\;[\mathfrak{c}^2_p]\;,\;\bullet\;]$ is linear functional on
$\det(H_p(C)^{\ast})\equiv \det(H_p(C))^{-1},$ where the exponent
$-1$ denotes the dual of the space. \\

This finishes the proof of Theorem~\ref{Torsion_Is_In_The_Dual}.
\end{proof}

\section{Symplectic Chain Complex}

\begin{defn}
$C_{\ast}:0\to C_n\stackrel{\partial_n}{\to} C_{n-1}\to\cdots\to
C_{\frac{n}{2}}\to\cdots \to C_1\stackrel{\partial_1}{\to} C_0\to 0$
is a \emph{symplectic chain complex}, if
\begin{itemize}
\item $n\equiv 2 ( \text{mod} 4 )$ and
\item there exist non-degenerate anti-symmetric $\partial-$compatible bilinear
maps i.e. $\omega_{p,n-p}:C_p\times C_{n-p}\to \mathbb{R}$ s.t.
$\omega_{p,n-p}(a,b)=(-1)^{p(n-p)}\omega_{n-p,p}(b,a)$ and
$\omega_{p,n-p}(\partial_{p+1}a,b)=(-1)^{p+1}\omega_{p+1,n-(p+1)}(a,\partial_{n-p}b).$
\end{itemize}
\end{defn}

In the definition, since $n\equiv 2 (\text{mod} 4)$ i.e. $n$ is even
and $\dfrac{n}{2}$ is odd,
$\omega_{p,n-p}(a,b)=(-1)^{p}\omega_{n-p,p}(b,a).$\\

Using the $\partial-$compatibility of the non-degenerate
anti-symmetric bilinear maps $\omega_{p,n-p}:C_p\times C_{n-p}\to
\mathbb{R},$ one can easily extend these to homologies. Namely,
\begin{lem}\label{Non-Degeneracy_Of [Omega]}
The bilinear map $[\omega_{p,n-p}]:H_p(C)\times H_{n-p}(C)\to
\mathbb{R}$ defined by
$[\omega_{p,n-p}]([x],[y])=\omega_{p,n-p}(x,y)$ is anti-symmetric and non-degenerate.\\
\end{lem}

\begin{proof}
For the well-definiteness, let $x,x'$ be in $\ker\partial_p$ with
$x-x'=\partial_{p+1}x''$ for some $x''\in C_{p+1}$ and let $y,y'$ be
in $\ker\partial_{n-p}$ with $y-y'=\partial_{n-p+1}y''$ for some
$y''\in C_{n-p+1}.$ Then from the bilinearity and
$\partial-$compatibility, $[\omega_{p,n-p}]([x],[y])$ is equal to
$\omega_{p,n-p}(x',y')+(-1)^p\omega_{p-1,n-p+1}(\partial_{p}x',y'')+
(-1)^{p+1}\omega_{p+1,n-p-1}(x'',\partial_{n-p}y')+
(-1)^{p+1}\omega_{p+1,n-p-1}(x'',\partial_{n-p}\circ\partial_{n-p+1}y'')
=\omega_{p,n-p}(x',y').$\\

Assume for some $[y_0]\in H_{n-p}(C),$
$[\omega_{p,n-p}]([x],[y_0])=0$ for all $[x]\in H_p(C).$\\

\begin{lem}\label{y_0_in_Image}
$y_0$ is in $\text{Im}\partial_{n-p+1}.$
\end{lem}

\begin{proof}
Let $\varphi:\dfrac{C_p}{Z_p}\to \mathbb{R}$ be defined by
$\varphi(x+Z_p)=\omega_{p,n-p}(x,y_0).$ This is a well-defined
linear map because if $x-x'\in Z_p,$ then
$\omega_{p,n-p}(x,y_0)-\omega_{p,n-p}(x',y_0)=[\omega_{p,n-p}]([x-x'],[y_0])$
equals to 0. By the $1^{\text{st}}$ isomorphism theorem,
$\dfrac{C_p}{Z_p}\stackrel{\overline{\partial_p}}{\cong}
\text{Im}\partial_p=B_{p-1},$ where $\overline{\partial_p}(x+Z_p)$
is $\partial_p(x).$\\

Consider the linear functional $\widetilde{\varphi}:=\varphi\circ
(\overline{\partial_p})^{-1}$ on $B_{p-1},$ where
$(\overline{\partial_p})^{-1}(\partial_p y)=y+Z_p.$ Extend
$\widetilde{\varphi}$ to $\widehat{\varphi}:C_{p-1}=B_{p-1}\oplus
\dfrac{C_{p-1}}{B_{p-1}}\to \mathbb{R}$  as zero on complement.
Since $\omega_{p-1,n-p+1}:C_{p-1}\times C_{n-p+1}\to\mathbb{R}$ is
non-degenerate, it induces an isomorphism between the dual space
$C^{\ast}_{p-1}$ of $C_{p-1}$ and $C_{n-p+1}.$ Therefore, there
exists a unique $u_0\in C_{n-p+1}$ such that
$\widehat{\varphi}(\cdot)=\omega_{p-1,n-p+1}(\cdot,u_0).$\\

For $x\in C_p,$ $v=\partial_p(x)$ is in $B_{p-1}$ Then, on one hand,
$\widehat{\varphi}(v)$ is $\omega_{p-1,n-p+1}(\partial_p x,u_0)$ or
$(-1)^p\omega_{p,n-p}( x,\partial_{n-p+1}u_0)$ by the
$\partial-$compatibility. On the other hand, by the construction of
$\widehat{\varphi},$ $\widehat{\varphi}(v)=\omega_{p,n-p}(x,y_0).$
So, $\omega_{p,n-p}(x,y_0)$ is $\omega_{p,n-p}(x,(-1)^p\partial_{n-p+1}u_0)$ for all $x\in C_p.$\\

The nondegeneracy of $\omega_{p,n-p}$ finishes the proof of
Lemma~\ref{y_0_in_Image}.
\end{proof}
This concludes the proof of Lemma~\ref{Non-Degeneracy_Of [Omega]}\\
\end{proof}

We will define $\omega-$compatibility for bases in a symplectic
chain complex.
\begin{defn}
Let $C_{\ast}:0\to C_n\stackrel{\partial_n}{\to} C_{n-1}\to\cdots\to
C_{\frac{n}{2}}\to\cdots C_1\stackrel{\partial_1}{\to} C_0\to 0$ be
a symplectic chain complex. Bases
$\mathfrak{o}_p,\mathfrak{o}_{n-p}$ of $C_p,C_{n-p}$ are
$\omega-$\emph{compatible} if the matrix of
$\omega_{p,n-p}$ in bases $\mathfrak{o}_p,\mathfrak{o}_{n-p}$ is\\
$$\left\{%
\begin{array}{ll}
    Id_{k\times k} & ; p\ne \frac{n}{2}  \\
    \left[%
\begin{array}{cc}
  O_{m\times m} & Id_{m\times m} \\
  -Id_{m\times m} & 0_{m\times m} \\
\end{array}%
\right] & ; p=\frac{n}{2} \\
\end{array}%
\right. $$ \\
where $k$ is $\dim C_p=\dim C_{n-p}$ and $2m=\dim C_{\frac{n}{2}}.$
In the same way, considering $[\omega_{p,n-p}]:H_p(C)\times
H_{n-p}(C)\to \mathbb{R},$ we can also define
$[\omega_{p,n-p}]-$compatibility of bases $\mathfrak{h}_p,$
$\mathfrak{h}_{n-p}$ of $H_p(C),$ $H_{n-p}(C).$
\end{defn}

In the next result, we will explain how a general symplectic chain
complex $C_{\ast}$ can be splitted $\omega-$orthogonally as a direct
sum of an exact and $\partial-$zero symplectic complexes.

\begin{thm}\label{Unnatural_Splitting_Of_Symplectic_Complex}
Let $C_{\ast}:0\to C_n\stackrel{\partial_n}{\to} C_{n-1}\to\cdots\to
C_1\stackrel{\partial_1}{\to} C_0\to 0$ be a symplectic chain
complex. Assume $\mathfrak{o}_p,\mathfrak{o}_{n-p}$
$\omega-$compatible. Then $C_{\ast}$ can be splitted as a direct sum
of symplctic complexes $C'_{\ast},C''_{\ast},$ where $C'_{\ast}$ is
exact, $C''_{\ast}$ is $\partial-$zero and $C'_{\ast}$ is
perpendicular to $C''_{\ast}.$
\end{thm}

\begin{proof}
Start with the following short-exact sequence
$$\begin{array}{cccccccc}
  0 & \to & \ker\partial_p & \hookrightarrow &
  C_p & \stackrel{\partial_p}{\rightarrow }  &
  \text{Im} \partial_p  & \to  0 \\
  0 & \to & \text{Im} \partial_{p+1} &
  \hookrightarrow & \ker\partial_p &
  \stackrel{\pi_p}{\rightarrow}  & H_p  & \to  0. \\
 \end{array}$$
Consider the section $\ell_p:\text{Im}\partial_p\to C_p$ defined by
$\ell_p(\partial_p x)=x$ for $\partial_p x\neq 0,$ and $s_p:H_p\to
\ker\partial_p$ by $s_p([x])=x.$\\

As $C_p$ disjoint union of $\text{Im}\partial_{p+1},$ $s_p(H_p),$
and $\ell_p(\text{Im}\partial_p),$ the basis $\mathfrak{o}_p$ of
$C_p$ has three blocks $\mathfrak{o}^1_p, \mathfrak{o}^2_p,
\mathfrak{o}^3_p,$ where $\mathfrak{o}^1_p$ is a basis for
$\text{Im}\partial_{p+1},$ $\mathfrak{o}^2_p$ generates $s_p(H_p)$
the rest  of $\ker\partial_p,$ i.e. $[\mathfrak{o}^2_p]$ generates
$H_p(C),$ and $\partial_p\mathfrak{o}^3_p$ is a basis for
$\text{Im}\partial_p.$ Similarly,
$\mathfrak{o}_{n-p}=\mathfrak{o}^1_{n-p}\sqcup
\mathfrak{o}^2_{n-p}\sqcup \mathfrak{o}^3_{n-p}.$ Because
$[\omega]_{p,n-p}:H_p(C)\times H_{n-p}(C)\to \mathbb{R}$ defined by
$[\omega]_{p,n-p}([a],[b])=\omega_{p,n-p}(a,b)$ is non-degenerate
and bases $\mathfrak{o}_p$ $\mathfrak{o}_{n-p}$ of $C_p,C_{n-p}$ are
$\omega$-compatible, $\omega_{p,n-p}(\cdot,s_{n-p}(H_{n-p})):C_p\to
\mathbb{R}$ vanishes on $\text{Im}\partial_{p+1}\oplus
\ell_p(\text{Im}\partial_p).$ Likewise,
$\omega_{p,n-p}(s_p(H_p(C)),\cdot):C_{n-p}\to \mathbb{R}$ vanishes
on $\text{Im}\partial_{n-p+1}\oplus
\ell_{n-p}(\text{Im}\partial_{n-p}).$\\

Set $C'_p=\text{Im}\partial_{p+1}\oplus \ell_p(\text{Im}\partial_p)$
and $C''_p=s_p(H_p).$ Note that $C'_p$ with basis
$\mathfrak{o}^1_p\sqcup \mathfrak{o}^3_p$ and $C''_{n-p}$ with basis
$\mathfrak{o}^2_{n-p}$ are $\omega-$orhogonal to each other. Hence,
$(C'_{\ast},\partial),$ $(C''_{\ast},\partial)$ will be the desired
splitting, where we consider the corresponding restrictions of
$\omega_{p,n-p}:C_p\times C_{n-p}\to \mathbb{R}.$\\

Clearly, $(C''_{\ast},\partial)$ is $\partial-$zero for $C''_p$
being subspace of $\ker\partial_p.$ Since
$[\omega_{p,n-p}]:H_p(C)\times H_{n-p}(C)\to \mathbb{R}$ is
non-degenerate, the restriction $\omega_{p,n-p}:C''_p\times
C''_{n-p}\to \mathbb{R}$ is also non-degenerate. Being the
restriction of $\omega_{p,n-p},$ it is also $\partial-$compatible.
Hence $C''_{\ast}$ becomes symplectic chain complex with
$\partial-$zero.

In the sequence
 $C'_{p+1}\stackrel{\partial_{p+1}}{\to}
C'_{p}\stackrel{\partial_{p}}{\to} C'_{p-1},$ first map
$\partial_{p+1}$ sends $\text{Im}\partial_{p+2},$
$\ell_{p+1}(\text{Im}\partial_{p+1})$  to zero and
$\text{Im}\partial_{p+1},$ respectively. Hence,
$\ker\{\partial_{p+1}:C'_{p+1}\to C'_p\}$ equals to
$\text{Im}\{\partial_{p+2}:C_{p+2}\to C_{p+1}\}$ and
$\text{Im}\{\partial_{p+1}:C'_{p+1}\to C'_p\}$ is
$\text{Im}\{\partial_{p+1}:C_{p+1}\to C_p\}.$ Similarly,
$\ker\{\partial_{p}:C'_{p}\to
C'_{p-1}\}=\text{Im}\{\partial_{p+1}:C_{p+1}\to C_{p}\}$ and
$\text{Im}\{\partial_{p}:C'_{p}\to C'_{p-1}\}$ is
$\text{Im}\{\partial_{p}:C_{p}\to C_{p-1}\}.$ Thus, $C'_{\ast}$ is
exact. \\

Moreover, since $\omega_{p,n-p}:C_p\times C_{n-p}\to \mathbb{R}$ is
non-degenerate, and $C'_p,C'_{n-p}$ are $\omega-$perpendicular to
$C''_{n-p},$ $C''_{p},$ respectively, $\omega_{p,n-p}:C'_p\times
C'_{n-p}\to \mathbb{R}$ is non-degenerate. Also, it is
$\partial-$compatible for being restriction of the
$\partial-$compatible map $\omega_{p,n-p}:C_p\times C_{n-p}\to
\mathbb{R}.$\\

This concludes the  proof of
Theorem~\ref{Unnatural_Splitting_Of_Symplectic_Complex}
\end{proof}

Above theorem is a special case of
Theorem~\ref{Unnatural_Splitting_Of_General_Complex}. The only only
difference is using $\omega-$compatible bases $\mathfrak{o}_p$ the
splitting is $\omega-$orthogonal, too.\\

We will now explain how the torsion of a symplectic complex with
$\partial-$zero is connected with Pfaffian of the anti-symmetric
$[\omega_{\frac{n}{2},\frac{n}{2}}]:H_{\frac{n}{2}}(C)\times
H_{\frac{n}{2}}(C)\to \mathbb{R}.$ Then, Pfaffian will be defined.
After that, we will give the relation for a general symplectic
complex.

\begin{thm}\label{Tor_of_Bondary_0_Symplectic}
Let $C_{\ast}$ be symplectic chain complex with $\partial-$zero. Let
$\mathfrak{h}_p$ be a basis for $H_p.$ Assume the bases
$\mathfrak{o}_p,\mathfrak{o}_{n-p}$ of $C_p,C_{n-p}$ are
$\omega-$compatible with the property that the bases
$\mathfrak{o}_{\frac{n}{2}}$ and $h_{\frac{n}{2}}$ of
$H_{\frac{n}{2}}(C)$ are in the same orientation class. Then,
$$\text{Tor}(C_{\ast},\{\mathfrak{o}_p\}_{p=0}^n,\{\mathfrak{h}_p\}_{p=0}^n)=
\left(\displaystyle\prod_{p=0}^{\frac{n}{2}-1} (\det
[\omega_{p,n-p}])^{(-1)^p}\right)\cdot
   \left(\sqrt{\det[\omega_{\frac{n}{2},\frac{n}{2}}]}\;\right)^{(-1)^{\frac{n}{2}}},$$
where $\det[\omega_{p,n-p}]$ is the determinant of the matrix of the
non-degenerate pairing $[\omega_{p,n-p}]:H_p(C)\times H_{n-p}(C)\to
\mathbb{R}$ in bases $\mathfrak{h}_p,\mathfrak{h}_{n-p}.$
\end{thm}

\begin{proof}
$C_{\ast}$ is $\partial-$zero complex, so all $B_p$'s are zero and
$Z_p=C_p.$ In particular, $H_p$ is equal to $C_p/\{0\}$ and hence
the basis $\mathfrak{h}_p$ of $H_p$ can also be considered as a
basis in $C_p.$ Recall
$\text{Tor}(C_{\ast},\{\mathfrak{o}_p\}_{p=0}^n,\{\mathfrak{h}_p\}_{p=0}^n)$
is defined as the alternating product
$$\prod_{p=0}^{n}[\mathfrak{o}_p,\mathfrak{h}_p]^{(-1)^p}=[\mathfrak{o}_0,\mathfrak{h}_0]^{(-1)^0}
\cdots
[\mathfrak{o}_{\frac{n}{2}},\mathfrak{h}_{\frac{n}{2}}]^{(-1)^{\frac{n}{2}}}
\cdots [\mathfrak{o}_n,\mathfrak{h}_n]^{(-1)^n},$$ of the
determinants $[\mathfrak{o}_p,\mathfrak{h}_p]$ of the
change-base-matrices
 from $\mathfrak{h}_p$ to
$\mathfrak{o}_p.$ If we combine the terms symmetric with the middle
term
$[\mathfrak{o}_{\frac{n}{2}},\mathfrak{h}_{\frac{n}{2}}]^{(-1)^{\frac{n}{2}}},$
torsion becomes
$$\displaystyle
\left(\prod_{p=0}^{\frac{n}{2}-1}[\mathfrak{o}_{p},\mathfrak{h}_{p}]^{(-1)^{p}}
[\mathfrak{o}_{n-p},\mathfrak{h}_{n-p}]^{(-1)^{n-p}}\right)
[\mathfrak{o}_{\frac{n}{2}},\mathfrak{h}_{\frac{n}{2}}]^{(-1)^{\frac{n}{2}}}.$$

Moreover, note that $[\mathfrak{o}_{p},\mathfrak{h}_{p}]^{(-1)^{p}}
[\mathfrak{o}_{n-p},\mathfrak{h}_{n-p}]^{(-1)^{n-p}}= \large\{\;
[\mathfrak{o}_{p},\mathfrak{h}_{p}]
[\mathfrak{o}_{n-p},\mathfrak{h}_{n-p}]\;\large\}^{(-1)^{p}}$ for
$n$ being even. Let
$T_{\mathfrak{h}_p}^{\mathfrak{o}_p},T_{\mathfrak{h}_{n-p}}^{\mathfrak{o}_{n-p}}$
denote the change-base-matrices from $\mathfrak{h}_p$ to
$\mathfrak{o}_p$ of $C_p,$ and from $\mathfrak{h}_{n-p}$ to
$\mathfrak{o}_{n-p}$ of $C_{n-p}$ respectively, i.e.
$h^i_p=\displaystyle\sum_{\alpha}(T_{\mathfrak{h}_p}^{\mathfrak{o}_p})_{\alpha
i}\;o^\alpha_p$ and
$h^j_{n-p}=\displaystyle\sum_{\beta}(T_{\mathfrak{h}_{n-p}}^{\mathfrak{o}_{n-p}})_{\beta
j}\;o^\beta_{n-p},$ where $h^i_p$ is the $i^{\text{th}}$-element of
the basis
$\mathfrak{h}_p.$\\

If $A$ and $B$ are the matrices of $\omega_{p,n-p}$ in the bases
$\mathfrak{h}_p,\mathfrak{h}_{n-p},$ and in the bases
$\mathfrak{o}_p,\mathfrak{o}_{n-p},$ respectively, then
$A=(T_{\mathfrak{h}_p}^{\mathfrak{o}_p})^{\text{transpose}}\;B\;
T_{\mathfrak{h}_{n-p}}^{\mathfrak{o}_{n-p}}.$ By the
$\omega-$compatibility of the bases
$\mathfrak{o}_p,\mathfrak{o}_{n-p},$ the matrix $B$ is equal to
$Id_{k\times k},$ $\left[%
\begin{array}{cc}
0_{m\times m} & Id_{m\times m} \\
-Id_{m\times m} & 0_{m\times m} \\
\end{array}%
\right]$  for $p\ne\frac{n}{2},$ $p=\frac{n}{2},$ respectively,
where $k$ is $\dim C_p=\dim C_{n-p}$ and $2m=\dim C_{\frac{n}{2}}.$
Clearly,
determinant of $B$ is $1^k=(-1)^m(-1)^m$ or $1.$\\

Hence, $\det A$ equals $\det T_{\mathfrak{h}_p}^{\mathfrak{o}_p}\det
T_{\mathfrak{h}_{n-p}}^{\mathfrak{o}_{n-p}}$ or
$[\mathfrak{o}_{p},\mathfrak{h}_{p}]
[\mathfrak{o}_{n-p},\mathfrak{h}_{n-p}]$ for all $p.$ In particular,
for $p=\frac{n}{2},$ it is
$[\mathfrak{o}_{\frac{n}{2}},\mathfrak{h}_{\frac{n}{2}}]^2.$ Since
$2m$ is even, and $\omega_{\frac{n}{2},\frac{n}{2}}$ is
non-degenerate skew-symmetric, the determinant $\det
A_{\frac{n}{2}}$ is positive actually equals to
$\text{Pfaf}(\omega_{\frac{n}{2},\frac{n}{2}})^2,$ and thus
$[\mathfrak{o}_{\frac{n}{2}},\mathfrak{h}_{\frac{n}{2}}]=\pm\sqrt{\det
A_{\frac{n}{2}}}.$ Because
$\mathfrak{o}_{\frac{n}{2}},\mathfrak{h}_{\frac{n}{2}}$ are in the
same orientation class, then
$[\mathfrak{o}_{\frac{n}{2}},\mathfrak{h}_{\frac{n}{2}}]=\sqrt{\det
A_{\frac{n}{2}}}.$\\

The proof is finished by the fact
$\omega_{p,n-p}(h^i_p,h^j_{n-p})=[\omega_{p,n-p}](h^i_p,h^j_{n-p}).$
\end{proof}

\begin{thm}\label{Tor_of_an Exact_Symplectic}
Let $C_{\ast}$ be an exact sypmlectic chain complex. If
$\mathfrak{c}_p,\mathfrak{c}_{n-p}$ are bases for $C_p,C_{n-p},$
respectively, then
$\text{Tor}(C_{\ast},\{\mathfrak{c}_p\}_{p=0}^n,\{0\}_{p=0}^n)=1.$
\end{thm}

\begin{proof}
From the exactness of $C_{\ast},$ we have $H_p=0$ or $\ker
\partial_p=\text{Im}\partial_{p+1}.$ Using the
short-exact sequence
$$0\to \ker\partial_p\hookrightarrow C_p\to
\text{Im}\partial_p\to 0,$$ we also have $C_p=\ker\partial_p\oplus
\ell_p(\text{Im}\partial_p),$
 where we consider the section
$\ell_p(\partial_p x)=x$ for $\partial_p x\ne0.$\\

Let $\mathfrak{o}_p,\mathfrak{o}_{n-p}$ be $\omega-$compatible bases
for $C_p,C_{n-p},$ respectively. We can break
$\mathfrak{o}_p=\mathfrak{o}^1_p\sqcup \mathfrak{o}^3_p,$ where
$\mathfrak{o}^1_p$ generates
$\ker\partial_p=\text{Im}\partial_{p+1},$ and
$\partial_p\mathfrak{o}^3_p$ generates $\text{Im}\partial_p.$
Similarly, $\mathfrak{o}_{n-p}=\mathfrak{o}^1_{n-p}\sqcup
\mathfrak{o}^3_{n-p},$ where $\mathfrak{o}^1_{n-p}$ generates
$\ker\partial_{n-p}=\text{Im}\partial_{n-p+1},$ and
$\partial_{n-p}\mathfrak{o}^3_{n-p}$ generates
$\text{Im}\partial_{n-p}.$ Since $\omega_{p,n-p}:C_p\times
C_{n-p}\to \mathbb{R}$ is  non-degenerate, $\partial-$compatible,
then $\omega_{p,n-p}(\mathfrak{o}^1_p,\mathfrak{o}^1_{n-p})=0,$ and
$\omega_{p,n-p}(\mathfrak{o}^1_p,\mathfrak{o}^3_{n-p})$ does not
vanish. Also by the $\omega-$compatibility of
$\mathfrak{o}_p,\mathfrak{o}_{n-p},$ for every $i$ there is unique
$j_i$ such that
$\omega_{p,n-p}((\mathfrak{o}^1_p)_i,(\mathfrak{o}^3_{n-p})_\alpha)=\delta_{j_i,\alpha}.$
Likewise, for every $k$ there is unique $q_k$ such that
$\omega_{p,n-p}((\mathfrak{o}^3_p)_k,(\mathfrak{o}^1_{n-p})_{\beta})=\delta_{q_k,\beta}.$\\

Recall torsion is independent of bases $\mathfrak{b}_p$ for
$\text{Im}\partial_p$ and section $\text{Im}\partial_p\to C_p.$ Let
$A_p$ be the determinant of the matrix of $\omega_{p,n-p}$ in bases
$\mathfrak{c}_p$ $\mathfrak{c}_{n-p},$ and let $O_p$ be the
determinant of the matrix of $\omega_{p,n-p}$ in bases
$\mathfrak{o}^1_p\sqcup \mathfrak{o}^3_p,$
$\mathfrak{o}^1_{n-p}\sqcup \mathfrak{o}^3_{n-p},$
 Since the set
$\partial_p\mathfrak{o}^3_p=\{\partial_p((\mathfrak{o}^3_p)_1),\cdots,
\partial_p((\mathfrak{o}^3_p)_{\alpha})\}$
generates $\text{Im}\partial_p,$ so does the set
$\{\partial_p(A_pO_p(\mathfrak{o}^3_p)_1),\partial_p((\mathfrak{o}^3_p)_2)\cdots,
\partial_p((\mathfrak{o}^3_p)_{\alpha})\}.$ Hence, image of the
latter set under $\ell_p,$ namely,
$\widetilde{\mathfrak{o}^3_p}=\{A_p\cdot
O_p\cdot(\mathfrak{o}^3_p)_1,(\mathfrak{o}^3_p)_2\cdots,
(\mathfrak{o}^3_p)_{\alpha}\}$
 will also be basis for $\ell_p(\text{Im}\partial_p).$
Keeping $\widetilde{\mathfrak{o}^3_{n-p}}$ as
$\mathfrak{o}^3_{n-p},$ we have
 $$\left[%
\begin{array}{cc}
  \omega_{p,n-p} &\text{in} \\
\\
  \mathfrak{o}^1_p\sqcup\widetilde{\mathfrak{o}^3_p}&\mathfrak{o}^1_{n-p}\sqcup\mathfrak{o}^3_{n-p} \\
\end{array}%
\right]=
(T_{\mathfrak{o}^1_p\sqcup\widetilde{\mathfrak{o}^3_p}}^{\mathfrak{c}_p})^{\text{transpose}}
\left[\begin{array}{cc}
  \omega_{p,n-p} &\text{in} \\
  \\
  \mathfrak{c}_p &\mathfrak{c}_{n-p} \\
\end{array}%
\right]
T_{\mathfrak{o}^1_{n-p}\sqcup\mathfrak{o}^3_{n-p}}^{\mathfrak{c}_{n-p}}.
 $$\\
Determinant of left-hand-side is $A_p\cdot O_p\cdot O_p,$ or $A_p$
because of the determinant of $\omega_{p,n-p}$ in the
$\omega-$compatible bases $\mathfrak{o}_p,\mathfrak{o}_{n-p}.$ Thus,
for $p\ne \frac{n}{2},$ we obtained that
$[\mathfrak{c}_p,\mathfrak{o}^1_p\sqcup\widetilde{\mathfrak{o}^3_p}]
[\mathfrak{c}_{n-p},\mathfrak{o}^1_{n-p}\sqcup\mathfrak{o}^3_{n-p}]=1.$
\\

For $p=\frac{n}{2},$ we can prove the same property as follows.
Since $\frac{n}{2}$ is odd,
$\omega_{\frac{n}{2},\frac{n}{2}}:C_{\frac{n}{2}}\times
C_{\frac{n}{2}}\to \mathbb{R}$ is non-degenerate and alternating,
then
 the matrix of $\omega_{\frac{n}{2},\frac{n}{2}}$ in any basis of $C_{\frac{n}{2}}$
will be an invertible $2m\times 2m$ skew-symmetric matrix $X$ with
real entries, where $2m=\dim C_{\frac{n}{2}}.$ Actually, we can find
an orthogonal $2m\times 2m$ matrix $Q$ with real entries so that

$$QXQ^{-1}=\text{diag}\left(\left(%
\begin{array}{cc}
  0 & a_1 \\
  -a_1 & 0 \\
\end{array}%
\right),\cdots,\left(%
\begin{array}{cc}
  0 & a_m \\
  -a_m & 0 \\
\end{array}%
\right)\right).$$\\

 So, the determinant of $\omega_{\frac{n}{2},\frac{n}{2}}$ in any basis will
be positive, in particular, the determinants $A_{\frac{n}{2}},
O_{\frac{n}{2}}$ of $\omega_{\frac{n}{2},\frac{n}{2}}$ in basis
$\mathfrak{c}_{\frac{n}{2}},$ $\mathfrak{o}^1_{\frac{n}{2}}\sqcup
\mathfrak{o}^3_{\frac{n}{2}}$ respectively will be positive. Having
noticed that, let $\widetilde{\mathfrak{o}^3_{\frac{n}{2}}}=\{
\sqrt{A_{\frac{n}{2}}}\cdot \sqrt{O_{\frac{n}{2}}} \;\cdot
(\mathfrak{o}^3_{\frac{n}{2}})_1,(\mathfrak{o}^3_{\frac{n}{2}})_2\cdots,
(\mathfrak{o}^3_{\frac{n}{2}})_{\alpha}\}.$\\

As explained above, on one side, we have that
$\det\left[%
\begin{array}{c}
  \omega_{\frac{n}{2},\frac{n}{2}}\; \text{in} \\
\mathfrak{o}^1_{\frac{n}{2}}\sqcup\widetilde{\mathfrak{o}^3_{\frac{n}{2}}}
  \end{array}%
\right]$ is equal to
$\sqrt{A_{\frac{n}{2}}}\;\cdot\sqrt{A_{\frac{n}{2}}}\;\sqrt{O_{\frac{n}{2}}}\;\cdot\sqrt{O_{\frac{n}{2}}}\;\det
\left[%
\begin{array}{cc}
  \omega_{\frac{n}{2},\frac{n}{2}} \;\text{in} \\
  \mathfrak{o}^1_{\frac{n}{2}}\sqcup\mathfrak{o}^3_{\frac{n}{2}} \\
\end{array}%
\right] $ or $A_{\frac{n}{2}}.$ On the other side, it is the product
$[\mathfrak{c}_{\frac{n}{2}},\mathfrak{o}^1_{\frac{n}{2}}\sqcup\widetilde{\mathfrak{o}^3_{\frac{n}{2}}}]\cdot
A_{\frac{n}{2}}\cdot
[\mathfrak{c}_{\frac{n}{2}},\mathfrak{o}^1_{\frac{n}{2}}\sqcup\widetilde{\mathfrak{o}^3_{\frac{n}{2}}}].$
Consequently, $[\mathfrak{c}_{\frac{n}{2}},\mathfrak{o}^1_{\frac{n}{2}}\sqcup\widetilde{\mathfrak{o}^3_{\frac{n}{2}}}]^2$ is equal to $1.$\\

If
$\mathfrak{o}^1_{\frac{n}{2}}\sqcup\widetilde{\mathfrak{o}^3_{\frac{n}{2}}}$
and $\mathfrak{c}_{\frac{n}{2}}$ are already in the same orientation
class, then
$[\mathfrak{c}_{\frac{n}{2}},\mathfrak{o}^1_{\frac{n}{2}}\sqcup\widetilde{\mathfrak{o}^3_{\frac{n}{2}}}]=1.$
If not, considering $\widetilde{\mathfrak{o}^3_{\frac{n}{2}}}$ as
$\{ -\sqrt{A_{\frac{n}{2}}}\cdot \sqrt{O_{\frac{n}{2}}} \;\cdot
(\mathfrak{o}^3_{\frac{n}{2}})_1,(\mathfrak{o}^3_{\frac{n}{2}})_2\cdots,
(\mathfrak{o}^3_{\frac{n}{2}})_{\alpha}\},$ we still have
$[\mathfrak{c}_{\frac{n}{2}},\mathfrak{o}^1_{\frac{n}{2}}\sqcup\widetilde{\mathfrak{o}^3_{\frac{n}{2}}}]=1.$\\

As a result, we proved that
\begin{eqnarray*}
 && \text{Tor}(C_{\ast},\{\mathfrak{c}_p\}_{p=0}^n,\{0\}_{p=0}^n) = \displaystyle
\prod_{p=0}^{n}[\mathfrak{c}_p,\mathfrak{o}^1_p\sqcup\widetilde{\mathfrak{o}^3_p}]
^{(-1)^p} \\
   &&= \prod_{p=0}^{\frac{n}{2}-1}\left([\mathfrak{c}_p,\mathfrak{o}^1_p\sqcup\widetilde{\mathfrak{o}^3_p}]
[\mathfrak{c}_{n-p},\mathfrak{o}^1_{n-p}\sqcup\mathfrak{o}^3_{n-p}]\right)^{(-1)^p}\cdot
[\mathfrak{c}_{\frac{n}{2}},\mathfrak{o}^1_{\frac{n}{2}}\sqcup\widetilde{\mathfrak{o}^3_{\frac{n}{2}}}]^{(-1)^{\frac{n}{2}}}= 1\\
  \end{eqnarray*}
\end{proof}

Before explaining the corresponding result for a general symplectic
complex, we would like to recall the Pfaffian of a skew-symmetric
matrix.\\

Let $V$ be an even dimensional vector space over reals. Let $\omega
:V\times V\to\mathbb{R}$ be a bilinear and anti-symmetric. If we fix
a basis for $V,$ then $\omega $ can be represented by a $2m\times
2m$ skew-symmetric matrix.\\

If $A$ any $2m\times 2m$ skew-symmetric matrix with real entries
then, by the spectral theorem of normal matrices, one can easily
find an orthogonal $2m\times 2m$-real matrix $Q$ so
that $QAQ^{-1}=\text{diag}\left(\left(%
\begin{array}{cc}
  0 & a_1 \\
  -a_1 & 0 \\
\end{array}%
\right),\cdots,\left(%
\begin{array}{cc}
  0 & a_m \\
  -a_m & 0 \\
\end{array}%
\right)\right),$ where $a_1,\cdots,a_m$ are positive real.
Thus, in particular, determinant of $A$ is non-negative.\\

\begin{defn} For $2m\times 2m$ real skew-symmetric
matrix $A,$ \emph{Pfafian} of $A$ will be $\sqrt{\det A}.$\\
\end{defn}

Actually, if $A=[a_{ij}]$ is any  $2m\times 2m$ skew-symmetric
matrix and if we let
$\omega_A=\sum_{i<j}a_{ij}\overrightarrow{e_i}\wedge
\overrightarrow{e_j},$ then we can also define $\text{Pfaf}(A)$ as
the coefficient of $\overrightarrow{e_1}\wedge\cdots\wedge
\overrightarrow{e_{2m}}$ in the product
$\dfrac{\overbrace{\omega_A\wedge\cdots\wedge\omega_A}^{m-times}}{m!}.$\\

For example, if $A$ is the matrix $\text{diag}\left(\left(%
\begin{array}{cc}
  0 & a_1 \\
  -a_1 & 0 \\
\end{array}%
\right),\cdots,\left(%
\begin{array}{cc}
  0 & a_m \\
  -a_m & 0 \\
\end{array}%
\right)\right),$ then $\omega_A$ is $\displaystyle\sum_{i=1}^m
a_i\cdot\overrightarrow{e_{2i-1}}\wedge \overrightarrow{e_{2i}}.$ An
easy computation shows that
$\underbrace{\omega_A\wedge\omega_A\wedge\cdots\wedge
\omega_A}_{m-times}$ equals to $m!\underbrace{(a_1\cdots
a_m)}_{\text{Pfaffian of A}}\overrightarrow{e_1}\wedge\cdots\wedge
\overrightarrow{e_{2m}}.$\\

For a general $2m\times 2m$ skew-symmetric $A,$ we can find an
orthogonal matrix $Q$ such that $QAQ^{-1}=\text{diag}\left(\left(%
\begin{array}{cc}
  0 & a_1 \\
  -a_1 & 0 \\
\end{array}%
\right),\cdots,\left(%
\begin{array}{cc}
  0 & a_m \\
  -a_m & 0 \\
\end{array}%
\right)\right).$ As a result,
$\underbrace{\omega_{QAQ^{-1}}\wedge\omega_{QAQ^{-1}}\wedge\cdots\wedge
\omega_{QAQ^{-1}}}_{m-times}$ equals to $m!\underbrace{(a_1\cdots
a_m)}_{\text{Pfaffian of }
QAQ^{-1}}\overrightarrow{e_1}\wedge\cdots\wedge
\overrightarrow{e_{2m}}$ i.e.
$\text{Pfaf}(QAQ^{-1})=\sqrt{\det(QAQ^{-1})}$ or $\sqrt{det(A)}.$\\

On the other hand, one can easily prove that for any $2m\times 2m$
skew-symmetric matrix $X$ and any $2m\times 2m$ matrix $Y,$
$\text{Pfaf}(YXY^t)$ is equal to $\text{Pfaf}(A)\det(B).$
Consequently, since $Q$ is orthogonal matrix, we can conclude that
$\text{Pfaf}(A)^2=\det(A)$ for any skew-symmetric $2m\times 2m$ real
matrix $A.$ In other words, both definitions coincide.\\

Using Pfafian, we can rephrase
Theorem~\ref{Tor_of_Bondary_0_Symplectic} as follows.\\

If $C_{\ast}$ is a symplectic chain complex with $\partial-$zero,
$\mathfrak{h}_p$ is a basis for $H_p,$
$\mathfrak{o}_p,\mathfrak{o}_{n-p}$ $\omega-$compatible bases for
$C_p,C_{n-p}$ so that $\mathfrak{h}_{\frac{n}{2}}$ and
$[\mathfrak{o}_{\frac{n}{2}}]$ are in the same orientation class,
then
$$\text{Tor}(C_{\ast},\{\mathfrak{o}_p\}_{p=0}^n,\{\mathfrak{h}_p\}_{p=0}^n)=
\left(\displaystyle\prod_{p=0}^{\frac{n}{2}-1} (\det
[\omega_{p,n-p}])^{(-1)^p}\right)\cdot
   \left(\text{Pfaf }[\omega_{\frac{n}{2},\frac{n}{2}}]\;\right)^{(-1)^{\frac{n}{2}}},$$
where $\text{Pfaf }[\omega_{\frac{n}{2},\frac{n}{2}}]$ is the
Pfafian of the matrix of the non-degenerate pairing
$[\omega_{\frac{n}{2},\frac{n}{2}}]:H_{\frac{n}{2}}(C)\times
H_{\frac{n}{2}}(C)\to \mathbb{R}$ in bases
$\mathfrak{h}_{\frac{n}{2}},\mathfrak{h}_{\frac{n}{2}}.$\\

\begin{thm}\label{Tor_Of_General_Symplectic}
For a general symplectic complex $C_{\ast},$ if $\mathfrak{c}_p,$
$\mathfrak{h}_p$ are bases for $C_p,$ $H_p,$ respectively, then
$$\text{Tor}(C_{\ast},\{\mathfrak{c}_p\}_{p=0}^n,\{\mathfrak{h}_p\}_{p=0}^n)=
\left(\displaystyle\prod_{p=0}^{\frac{n}{2}-1} (\det
[\omega_{p,n-p}])^{(-1)^p}\right)\cdot
   \left(\sqrt{\det[\omega_{\frac{n}{2},\frac{n}{2}}]}\;\right)^{(-1)^{\frac{n}{2}}},$$
where $\det[\omega_{p,n-p}]$ is the determinant of the matrix of the
non-degenerate pairing $[\omega_{p,n-p}]:H_p(C)\times H_{n-p}(C)\to
\mathbb{R}$ in bases $\mathfrak{h}_p,\mathfrak{h}_{n-p}.$\\
\end{thm}

\begin{proof}
Since $C_p$ is disjoint union $\text{Im}\partial_{p+1}\sqcup
s_p(H_p)\sqcup\ell_p(\text{Im}\partial_p),$ any basis
$\mathfrak{a}_p$ of $C_p$ has three parts
$\mathfrak{a}_p^1,\mathfrak{a}_p^2,\mathfrak{a}_p^3,$ where
$\mathfrak{a}_p^1$ is basis for $\text{Im}\partial_{p+1},$
$\mathfrak{a}_p^2$ generates $s_p(H_p)$ the rest of $\ker\partial_p$
i.e. $[\mathfrak{a}_p^2]$ generates $H_p(C),$ and
$\partial_p\mathfrak{a}_p^3$ is basis for $\text{Im}\partial_p,$
where $\ell_p:\text{Im}\partial_p\to C_p$ is the section defined by
$\ell_p(\partial_p x)=x$ for $\partial_p x\neq 0,$ and $s_p:H_p\to
\ker\partial_p$ by $s_p([x])=x.$\\

 If $\mathfrak{o}_p,$ $\mathfrak{o}_{n-p}$ are
 $\omega-$compatible bases for $C_p$ and $C_{n-p},$ then we can also
 write $\mathfrak{o}_p=\mathfrak{o}_p^1\sqcup\mathfrak{o}_p^2\sqcup\mathfrak{o}_p^3$
and
$\mathfrak{o}_{n-p}=\mathfrak{o}_{n-p}^1\sqcup\mathfrak{o}_{n-p}^2\sqcup\mathfrak{o}_{n-p}^3.$
We may assume  $[\mathfrak{o}_{\frac{n}{2}}]$ and
$\mathfrak{h}_{\frac{n}{2}}$ are in the same orientation class.
Otherwise, switch, say, the first element
$(\mathfrak{o}_{\frac{n}{2}})^1$ and the corresponding
$\omega-$compatible element $(\mathfrak{o}_{\frac{n}{2}})^{m+1}$
i.e.
$\omega_{\frac{n}{2},\frac{n}{2}}((\mathfrak{o}_{\frac{n}{2}})^1,(\mathfrak{o}_{\frac{n}{2}})^{m+1})=1,$
where $2m=\dim H_{\frac{n}{2}}(C).$ In this way, we still have
$\omega-$compatibility and moreover we can guarantee that
$[\mathfrak{o}_{\frac{n}{2}}],$  $\mathfrak{h}_{\frac{n}{2}}$ are in
the same
orientation class.\\

Using these $\omega-$compatible bases $\mathfrak{o}_p,$ as in
Theorem~\ref{Unnatural_Splitting_Of_Symplectic_Complex}, we have the
$\omega-$orthogonal splitting $C_{\ast}=C'_{\ast}\oplus C''_{\ast},$
where $C'_p$ and $C''_p$ are $\text{Im}(\partial_{p+1})\oplus
\ell_p(\text{Im}\partial_p),$ $s_p(H_p),$ and
$\ell_p:\text{Im}\partial_p\to C_p$ is the section defined by
$\ell_p(\partial_p x)=x$ for $\partial_p x\neq 0,$ and $s_p:H_p\to
\ker\partial_p$ by $s_p([x])=x.$\\

$C_p$ is the disjoint union $\text{Im}\partial_{p+1}\sqcup
s_p(H_p)\sqcup\ell_p(\text{Im}\partial_p),$ so the basis
$\mathfrak{c}_p$ of $C_p$ has also three blocks $\mathfrak{c}^1_p,
\mathfrak{c}^2_p, \mathfrak{c}^3_p,$ where $\mathfrak{c}^1_p$ is a
basis for $\text{Im}\partial_{p+1},$ $\mathfrak{c}^2_p$ generates
$s_p(H_p)$ the rest  of $\ker\partial_p,$ i.e. $[\mathfrak{c}^2_p]$
generates $H_p(C),$ and $\partial_p\mathfrak{c}^3_p$ is a basis for
$\text{Im}\partial_p.$\\

Consider the $\partial-$zero symplectic $C''_{\ast}$ with the
$\omega-$compatible bases $\mathfrak{o}^2_p,\mathfrak{o}^2_{n-p}.$
Note that by the $\partial-$zero property of $C''_{\ast},$
$H_p(C'')$ is $C''_p/0$ or $s_p(H_p(C)).$ Hence
$s_p(\mathfrak{h}_p)$ will be a basis $H_p(C'').$ Recall also that
$[\mathfrak{o}^2_{\frac{n}{2}}]$ and $\mathfrak{h}^2_{\frac{n}{2}}$
are in the same orientation class. Therefore, by
Theorem~\ref{Tor_of_Bondary_0_Symplectic}, we can conclude that
$$\text{Tor}(C''_{\ast},\{\mathfrak{o}^2_p\}_{p=0}^n,\{s_p(\mathfrak{h}_p)\}_{p=0}^n)=
\left(\displaystyle\prod_{p=0}^{\frac{n}{2}-1} (\det
[\omega_{p,n-p}])^{(-1)^p}\right)\cdot
   \left(\sqrt{\det[\omega_{\frac{n}{2},\frac{n}{2}}]}\;\right)^{(-1)^{\frac{n}{2}}},$$
where $\det[\omega_{p,n-p}]$ is the determinant of the matrix of the
non-degenerate pairing $[\omega_{p,n-p}]:H_p(C)\times H_{n-p}(C)\to
\mathbb{R}$ in bases $\mathfrak{h}_p,\mathfrak{h}_{n-p}.$\\

On the other hand, if $\mathfrak{c}'_p$  is any basis for $C'_p,$
then by Theorem~\ref{Tor_of_an Exact_Symplectic} the torsion
$\text{Tor}(C'_{\ast},\{c'_p\}_{p=0}^n,\{0\}_{p=0}^n)$ of the exact
symplectic complex $C'_{\ast}$
 is equal to
$1$ .\\

Let $A_p$ be the determinant of the change-base-matrix from
$\mathfrak{o}^2_p$ to $\mathfrak{c}^2_p.$ If we consider the basis
$\mathfrak{c}^1_p\sqcup (\dfrac{1}{A_p}\mathfrak{c}^3_p)$ for the
$C'_p,$ then for the short-exact sequence
$$0\to C''_{\ast}\to C_{\ast}=C'_{\ast}\oplus C''_{\ast}\to C'_{\ast}\to 0 $$
the bases $\mathfrak{o}^2_p,\mathfrak{c}_p,\mathfrak{c}^1_p\sqcup
(\frac{1}{A_p}\mathfrak{c}^3_p)$ of $C''_p,C_p,C'_p$ respectively
will be compatible i.e. the determinant of the change-base-matrix
from basis $\mathfrak{c}^1_{p}\sqcup \mathfrak{o}^2_{p}\sqcup
(\frac{1}{A_p}\mathfrak{c}^3_p)$ to
$\mathfrak{c}_p=\mathfrak{c}^1_{p}\sqcup \mathfrak{c}^2_{p}\sqcup
\mathfrak{c}^3_{p}$ is $1.$\\

 Thus, by Milnor's result Theorem~\ref{Mil},
$\text{Tor}(C_{\ast},\{\mathfrak{c}_p\}_{p=0}^n,\{\mathfrak{h}_p\}_{p=0}^n)$
is equal to the product of
$\text{Tor}(C''_{\ast},\{\mathfrak{o}^2_p\}_{p=0}^n,\{s_p(\mathfrak{h}_p)\}_{p=0}^n),$
$\text{Tor}(C'_{\ast},\{\mathfrak{c}^1_p\sqcup
(\frac{1}{A_p}\mathfrak{c}^3_p)\}_{p=0}^n,\{0\}_{p=0}^n),$ and
$\text{Tor}(\mathcal{H}_{\ast},\{s_p(\mathfrak{h}_p),\mathfrak{h}_p,0\}_{p=0}^n,\{0\}_{p=0}^{3n+2}),$
where $\mathcal{H}_{\ast}$ is the long-exact sequence  $0\to
H_n(C'')\to H_n(C)\to H_n(C')\to H_{n-1}(C'')\to \cdots\to
H_0(C'')\to H_0(C)\to H_0(C')\to 0$ obtained from the short-exact
sequence of complexes. Since $C'_{\ast}$ is exact, $H_p(C')$ are all
zero. So, using the isomorphisms $H_p(C)\to H_p(C'')=C''_p/0,$
namely $s_p$ as section, we can conclude that
$\text{Tor}(\mathcal{H}_{\ast},\{s_p(\mathfrak{h}_p),\mathfrak{h}_p,0\}_{p=0}^n,\{0\}_{p=0}^{3n+2})=1.$
From Theorem~\ref{Tor_of_an Exact_Symplectic}, we also obtain
$\text{Tor}(C'_{\ast},\{\mathfrak{c}^1_p\sqcup
(\dfrac{1}{A_p}\mathfrak{c}^3_p)\}_{p=0}^n,\{0\}_{p=0}^n)=1.$ \\

Therefore, we verified that
$$\text{Tor}(C_{\ast},\{\mathfrak{c}_p\}_{p=0}^n,\{\mathfrak{h}_p\}_{p=0}^n)=
\text{Tor}(C''_{\ast},\{\mathfrak{o}^2_p\}_{p=0}^n,\{s_p(\mathfrak{h}_p)\}_{p=0}^n).$$

This finishes the proof of Theorem~\ref{Tor_Of_General_Symplectic}.\\
\end{proof}

\providecommand{\bysame}{\leavevmode\hbox
to3em{\hrulefill}\thinspace}
\providecommand{\MR}{\relax\ifhmode\unskip\space\fi MR }
\providecommand{\MRhref}[2]{%
  \href{http://www.ams.org/mathscinet-getitem?mr=#1}{#2}
} \providecommand{\href}[2]{#2}

\end{document}